\documentclass[10pt,reqno]{amsart}
\usepackage{bbm}
\usepackage{mathrsfs}
\usepackage{amsfonts} 
\usepackage[dvipsnames,usenames]{color}
\usepackage{graphicx,latexsym,bm,amsmath,amssymb,verbatim,multicol,lscape}
\usepackage{pgfplots}
\theoremstyle{definition}

\theoremstyle{remark}

\numberwithin{equation}{section}

\begin{document}
\title[On deep holes of generalized projective Reed Solomon codes]
{On deep holes of generalized projective Reed-Solomon codes}
\author{Xiaofan Xu}
\address{Mathematical College, Sichuan University, Chengdu 610064, P.R. China}
\email{xxfscu@163.com}  
\author{Shaofang Hong}
\address{Mathematical College, Sichuan University, Chengdu 610064, P.R. China}
\email{sfhong@scu.edu.cn, s-f.hong@tom.com, hongsf02@yahoo.com}
\author{Yongchao Xu}
\address{Mathematical College, Sichuan University, Chengdu 610064, P.R. China}
\email{xyongchao@hotmail.com}
\thanks{S. Hong is the corresponding author and was supported partially by
National Science Foundation of China Grant \# 11371260.}

\keywords{generalized projective Reed-Solomon code, MDS code, deep hole,
Lagrange interpolation polynomial, generator matrix}
\subjclass[2000]{Primary 11B25, 11N13, 11A05}
\date{\today}%
\begin{abstract}
Determining deep holes is an important topic in decoding Reed-Solomon
codes. Cheng and Murray, Li and Wan, Wu and Hong investigated the error
distance of generalized Reed-Solomon codes. Recently, Zhang and Wan
explored the deep holes of projective Reed-Solomon codes. Let $l\ge 1$
be an integer and $a_1,\ldots,a_l$ be arbitrarily given $l$ distinct
elements of the finite field ${\bf F}_q$ of $q$ elements with
the odd prime number $p$ as its characteristic.
Let $D={\bf F}_q\backslash\{a_1,\ldots,a_l\}$
and $k$ be an integer such that $2\le k\le q-l-1$. In this paper,
we study the deep holes of generalized projective Reed-Solomon code
${\rm GPRS}_q(D, k)$ of length $q-l+1$ and dimension $k$ over ${\bf F}_q$.
For any $f(x)\in {\bf F}_q[x]$, we let $f(D)=(f(y_1),\ldots,f(y_{q-l}))$
if $D=\{y_1, ..., y_{q-l}\}$ and $c_{k-1}(f(x))$ be the coefficient
of $x^{k-1}$ of $f(x)$. By using D\"ur's theorem on the relation between
the covering radius and minimum distance of ${\rm GPRS}_q(D, k)$,
we show that if $u(x)\in {\bf F}_q[x]$ with
$\deg (u(x))=k$, then the received codeword $(u(D), c_{k-1}(u(x)))$
is a deep hole of ${\rm GPRS}_q(D, k)$ if and only if the sum
$\sum\limits_{y\in I}y$ is nonzero for any subset $I\subseteq D$ with
$\#(I)=k$. We show also that if $j$ is an integer with $1\leq j\leq l$
and
$u_j(x):= \lambda_j(x-a_j)^{q-2}+\nu_j x^{k-1}+f_{\leq k-2}^{(j)}(x)$
with $\lambda_j\in {\bf F}_q^*$, $\nu_j\in {\bf F}_q$ and
$f_{\leq{k-2}}^{(j)}(x)\in{\bf F}_q[x]$ being a polynomial
of degree at most $k-2$, then $(u_j(D), c_{k-1}(u_j(x)))$
is a deep hole of ${\rm GPRS}_q(D, k)$ if and only if the sum
$\binom{q-2}{k-1}(-a_j)^{q-1-k}\prod\limits_{y\in I}(a_j-y)+e$
is nonzero for any subset $I\subseteq D$ with $\#(I)=k$, where $e$
is the identity of the group ${\bf F}_q^*$. This implies that
$(u_j(D), c_{k-1}(u_j(x)))$ is a deep hole of ${\rm GPRS}_q(D, k)$
if $p|k$. We also deduce that $(u({\bf F}_q^*),\delta)$ is a deep
hole of the primitive projective Reed-Solomon code
${\rm PPRS}_q({\bf F}_q^*, k)$
if $u(x)=\lambda x^{q-2}+\delta x^{k-1}+f_{\leq{k-2}}(x)$ with
$\lambda\in {\bf F}_q^*$ and $\delta\in {\bf F}_q$. But
$(u({\bf F}_q^*), c_{k-1}(u(x)))$ is not a deep hole of
${\rm PPRS}_q({\bf F}_q^*, k)$ if $\deg (u(x))=k$.
\end{abstract}

\maketitle

\section{\bf Introduction and the statements of the main results}

Let ${\bf F}_q$ be the finite field of $q$ elements with $p$
as its characteristic. Let $n$ and $k$ be positive integers
such that $k<n$. Let $D=\{x_1,\cdots,x_n\}$ be a subset of
${\bf F}_q$, which is called the {\it evaluation set}.
The {\it generalized Reed-Solomon code} ${\rm GRS}_q(D, k)$
of length $n$ and dimension $k$ over ${\bf F}_q$  is defined by:
$${\rm GRS}_q(D, k): = \{(f(x_1), \ldots, f(x_n))
\in {\bf F}_q^{n} \mid f (x) \in{\bf F}_q[x], \deg f(x) \leq {k-1}\}.$$
Moreover, the {\it generalized projective Reed-Solomon code}
${\rm GPRS}_q(D, k)$ of length $n+1$ and dimension
$k$ over ${\bf F}_q$  is defined as follows:
$${\rm GPRS}_q(D, k): = \{(f(x_1), \ldots, f(x_n),c_{k-1}(f(x)))
\in {\bf F}_q^{n+1} \mid f (x) \in{\bf F}_q[x], \deg f(x) \leq {k-1}\},$$
where $c_{k-1}(f(x))$ is the coefficient of $x^{k-1}$ of $f(x)$.
If $D = {\bf F}_q^*$ , then it is called
{\it primitive projective Reed-Solomon code}, namely,
$${\rm PPRS}_q({\bf F}_q^*, k) := \{(f(1),\ldots, f(\alpha^{q-2}),c_{k-1}(f(x)))
\in {\bf F}_q^q\mid f (x) \in{\bf F}_q[x], \deg f(x) \leq {k-1}\},$$
where $\alpha$ is a primitive element of ${\bf F}_q$. If $D = {\bf F}_q$,
then it is called the {\it extended projective Reed Solomon code}.
For $u=(u_1,\ldots,u_n)\in{\bf F}_q^n$, $v=(v_1, \ldots, v_n)
\in{\bf F}_q^n$, the {\it Hamming distance} $d(u,v)$ is defined by
$$d(u,v):=\#\{1\le i\le n \mid u_i\neq{v_i},u_i\in{\bf F}_q,v_i\in{\bf F}_q\}.$$
For any $[n, k]_q$ linear code $C$, the {\it minimum distance} $d(C)$ is defined by
$$ d(C):= \min\{d(x, y)\mid x\in C, y\in C, x\neq y \},$$
where $d(\cdot , \cdot)$ denotes the {\it Hamming distance} of two words.
A linear $[n, k, d]$ code is called {\it maximum distance separable} (MDS)
code if $d=n-k+1$. The {\it error distance} to code $C$ of a received word
$u \in {\bf F}_q^n$ is defined by $$d(u, C):= \min_{v\in C}\{d(u, v)\}.$$
Clearly, $d(u, C) = 0$ if and only if $u \in C$.
The {\it covering radius} to code $C$ of a received word $u \in {\bf F}_q^n$
is defined by $$\rho(C):= \max\{d(u, C)\mid u\in{\bf F}_q^n\}.$$

The most important algorithmic problem in coding theory is
the maximum likelihood decoding (MLD):
Given a received word $u\in{\bf F}_q^n$, find a word $v \in C$
such that $d(u, v) = d(u, C)$, then we decode $u$ to $v$ \cite{[LW2]}.
Therefore, it is very crucial to decide $d(u, C)$ for the word $u$.
Guruswami and Sudan \cite{[GS]} provided a polynomial time list
decoding algorithm for the decoding of $u$ when $d(u,C)\le {n-\sqrt{nk}}$.
When the error distance increases, Guruswami and Vardy \cite{[GV]}
showed that that maximum-likelihood decoding is NP-hard for the
family of Reed-Solomon codes. We also notice that D\"ur \cite{[D]}
studied the Cauchy codes. In particular, D\"ur \cite{[D]} got the
relation between the covering radius and minimum distance of
${\rm GPRS}_q(D, k)$. When decoding the generalized projective Reed-Solomon
code $C$,for a received word $u = (u_1,\ldots, u_n,u_{n+1})\in {\bf F}_q^{n+1}$,
we define the {\it Lagrange interpolation polynomial} $u(x)$ of $u$ by
$$u(x):=\sum\limits _{i=1}^{n}u_i\prod\limits_{j=1 \atop j\neq i}^{n}
\frac{x-x_j}{x_i-x_j}\in {\bf F}_q[x], $$
i.e., $u(x)$ is the unique polynomial of degree $\deg u(x) \le n-1$
such that $u(x_i) = u_i$ for $1\leq i\leq n$ and $c_{k-1}(u(x)) = u_{n+1}$.
It is clear that $u\in C$ if and only if $d(u, C)=0$
if and only if $\deg u(x)\leq  {k-1}$. Equivalently, $u\not \in C$ if and
only if $d(u, C)\ge 1$ if and only if $k\le \deg u(x)\le n-1$.
Evidently, we have the following simple bounds of $d(u,{\rm GRS}_q(D, k))$
which are due to Li and Wan.\\

\noindent{\bf Theorem 1.1.} \cite{[LW1]} {\it Let $u$ be a received word
such that $u\not\in {\rm GRS}_q(D, k)$. Then}
$$n-\deg u(x) \leq d(u, {\rm GRS}_q(D, k))\leq n-k= \rho({\rm GRS}_q(D, k)).$$\\

Let $u \in {\bf F}_q^n$. If $d(u, {\rm GRS}_q(D, k)) = \rho({\rm GRS}_q(D, k))$,
then the received word $u$ is called a {\it deep hole} of ${\rm GRS}_q(D, k)$.
In 2007, Cheng and Murray \cite{[CM]} conjectured that
a word $u$ is a deep hole of ${\rm GRS}_q({\bf F}_q,k)$ if and only if
$u(x) = ax^k+f_{\leq{k-1}}(x)$, where $u(x)$ is
the Lagrange interpolation polynomial of the received word $u$
and $a\in{\bf F}_q^*$, $f_{\leq{k-1}}(x)\in{\bf F}_q[x]$
a polynomial of degree at most $k-1$. In 2012, Wu and Hong \cite{[WH]}
disproved this conjecture by giving a new class of deep holes for
Reed-Solomon codes ${\rm GRS}_q({\bf F}_q^*, k)$. In fact,
if $q\geq 4$ and $2\leq k\leq q-2$, then they showed
that the received word $u$ is a deep hole if its Lagrange
interpolation polynomial is $ ax^{q-2} +f_{\leq{k-1}}(x) $.
In \cite{[HW]}, Hong and Wu proved that the received word
$u$ is a deep hole of the generalized Reed-Solomon codes
${\rm GRS}_q(D,k)$ if its Lagrange interpolation polynomial
is $\lambda(x-a_i)^{q-2}+f_{\leq{k-1}}(x)$,
where $\lambda\in {\bf F}_q^*, a_i\in{\bf F}_q\backslash D$
and $f_{\leq{k-1}}(x)\in{\bf F}_q[x]$ a polynomial of
degree at most $k-1.$

Throughout this paper, we always let $l$ be a positive integer and
$a_1,\ldots,a_l$ be any fixed $l$ distinct elements of ${\bf F}_q$.
Let
$$D:={\bf F}_q\backslash\{a_1,\ldots,a_l\}.$$
We write
$$D:=\{y_1,\ldots,y_{q-l}\},$$
and for any $f(x)\in {\bf F}_q[x]$, we let
$$f(D):=(f(y_1),\ldots,f(y_{q-l})),$$
and use $c_{k-1}(f(x))$ to denote the coefficient of $x^{k-1}$ of
$f(x)$. Then we can rewrite the generalized projective
Reed-Solomon code ${\rm GPRS}(D, k)$ with evaluation set $D$ as
$${\rm GPRS}_q(D, k):=\{(f(D),c_{k-1}(f(x)))
\in {\bf F}_q^{q-l+1} \mid f (x) \in{\bf F}_q[x], \deg f(x) \leq {k-1}\}.$$

Let $u\notin {\rm GPRS}_q(D, k)$. If
$d(u,{\rm GPRS}_q(D, k))=\rho({\rm GPRS}_q(D, k))$,
then $u$ is also called a {\it deep hole} of generalized projective
Reed-Solomon code ${\rm GPRS}_q(D, k)$. In 2016, Zhang and
Wan \cite{[ZW]} studied the deep holes of projective Reed-Solomon
code ${\rm GPRS}({\bf F}_q, k)$. In fact, under the assumption
that the only deep holes of ${\rm GRS}_q({\bf F}_q, k)$ are
those received codewords whose Lagrange interpolation polynomials
are of degree $k$, they proved the following results
by solving a subset sum problem.\\

\noindent{\bf Theorem 1.2.} \cite{[ZW]}
{\it Let $q$ be an odd prime power. Assume that $3\leq k+1\leq p$
or $3\leq q-p+1\leq k+1\leq {q-2}$. Then the received
codeword $(f({\bf F}_q), c_{k-1}(f(x)))$ with $\deg f(x)=k$
is a deep hole of ${\rm GPRS}({\bf F}_q, k)$.}\\

\noindent{\bf  Theorem 1.3.} \cite{[ZW]}
{\it Let $\deg f(x)\geq k+1$ and $s:=\deg f(x)-k+1$.
If there are positive constants $c_1$ and $c_2$
such that $s<c_1\sqrt{q}, (\frac{s}{2}+2)\log_2(q)<k<c_2q$,
then $(f({\bf F}_q), c_{k-1}(f(x)))$
is not a deep hole of ${\rm GPRS}({\bf F}_q, k)$.}\\

In this paper, our main goal is to investigate the deep holes of the
generalized projective Reed-Solomon codes ${\rm GPRS}_q(D,k)$. Actually, we
will present characterizations for the received codewords of degrees $k$ and $q-2$
to be deep holes of generalized projective Reed-Solomon code ${\rm GPRS}_q(D,k)$.
The main results of this paper can be stated as follows.\\

\noindent{\bf Theorem 1.4.} {\it Let $q$ be a prime
power and $k$ and $l$ be positive integers
such that $q\geq 5$ and $2\leq k \leq \min(q-3, q-l-1)$.
Let $u(x)\in {\bf F}_q[x]$ with $\deg (u(x))=k$.
Then the received codeword $(u(D), c_{k-1}(u(x)))$
is a deep hole of the generalized projective Reed-Solomon code
${\rm {\rm GPRS}}_q(D, k)$ if and only if the sum $\sum\limits_{y\in I}y$
is nonzero for any subset $I\subseteq D$ with $\#(I)=k$.}\\

\noindent{\bf Theorem 1.5.} {\it Let $q$ be a prime power and
$k$ and $l$ be positive integers such that $q\geq 4$ and
$2\leq k \leq {q-l-1}$. Let $j$ be an integer with $1\leq j\leq l$
and let
$u_j(x):= \lambda_j(x-a_j)^{q-2}+\nu_j x^{k-1}+f_{\leq k-2}^{(j)}(x)$
with $\lambda_j\in {\bf F}_q^*$,
$\nu_j\in {\bf F}_q$ and $f_{\leq{k-2}}^{(j)}(x)\in{\bf F}_q[x]$
being a polynomial of degree at most $k-2$.
Then the received codeword $(u_j(D), c_{k-1}(u_j(x)))$
is a deep hole of the generalized projective
Reed-Solomon code ${\rm GPRS}_q(D, k)$ if and only if
the sum $\binom{q-2}{k-1}a_j^{q-1-k}\prod\limits_{y\in I}(y-a_j)+e$
is nonzero for any subset $I\subseteq D$ with $\#(I)=k$, where $e$ is
the identity of the multiplicative group ${\bf F}_q^*$.

Further, if $k\equiv 0\pmod p$, then the received codeword
$(u_j(D), c_{k-1}(u_j(x)))$ is a deep hole of ${\rm GPRS}_q(D, k)$.}\\

From Theorems 1.4 and 1.5, we can derive the following results on the
deep holes of the primitive projective Reed-Solomon codes. Note that
the proof of Theorem 1.6 relies also on a result on the zero subsets sum
of the group ${\bf F}_q^*$ (see Lemma 2.8 below). \\

\noindent{\bf Theorem 1.6.} {\it Let $q$ be an odd prime power
such that $q\geq 5$ and $2\leq k \leq {q-3}$. If
$u(x)=\lambda x^k+\gamma x^{k-1}+f_{\leq{k-2}}(x)$ with
$\lambda\in {\bf F}_q^*,\gamma\in {\bf F}_q $
and $f_{\leq{k-2}}(x)\in{\bf F}_q[x]$
being a polynomial of degree at most $k-2$,
then the received codeword $(u({\bf F}_q^*),\gamma)$
is not a deep hole of the primitive projective
Reed-Solomon code ${\rm PPRS}_q({\bf F}_q^*, k)$.}\\

\noindent{\bf Theorem 1.7.} {\it Let $q\geq 4$ and $2\leq k \leq {q-2}$. If
$u(x)=\lambda x^{q-2}+\delta x^{k-1}+f_{\leq{k-2}}(x)$ with
$\lambda\in {\bf F}_q^*, \delta\in {\bf F}_q$ and
$f_{\leq{k-2}}(x)\in{\bf F}_q[x]$
being a polynomial of degree at most $k-2$, then the received codeword
$(u({\bf F}_q^*),\delta)$ is a deep hole of the primitive projective
Reed-Solomon code ${\rm PPRS}_q({\bf F}_q^*, k)$.}\\

In the proofs of Theorems 1.4 and 1.5, the basic tools are
the MDS code and Vandemonde determinant. But we would also
like to point out that a key ingredient in the proofs is
D\"ur's theorem on the relation between the covering radius
and minimum distance of the generalized projective Reed-Solomon
code ${\rm GPRS}_q(D, k)$ (see Lemma 2.6 below). Another
important ingredient is a new result on the zero-sum problem
in the finite field that we will prove in the next section.

This paper is organized as follows. First of all, in Section 2,
we recall and prove several preliminary lemmas that are needed
in the proof of Theorems 1.4 and 1.5. Consequently,  in Section 3,
we use the lemmas presented in Section 2 to give the proofs of
Theorems 1.4 and 1.6. Finally, by using the results given
in Section 2, we supply in Section 4 the proofs of
Theorems 1.5 and 1.7.

\section{\bf Preliminary lemmas}

In this section, our main aim is to prove several lemmas that are needed in
the proof of Theorems 1.4 and 1.5. We begin with the following result on MDS
codes. \\

\noindent{\bf Lemma 2.1.} {\it Let $C$ be a MDS code and
$u_0\in C$ be a given codeword. Then the received codeword $u$ is
a deep hole of $C$ if and only if the received codeword $u+u_0$ is
a deep hole of $C$.}\\

{\it Proof.} First of all, let $u$ be a received codeword.
Then by the definition of deep hole, one knows
that $u$ is a deep hole of $C$ if and only if
$d(u, C)=\rho (C)$ with $\rho (C)$ being the covering radius of $C$,
if and only if
\begin{align}
\min_{v\in C}\{d(u, v)\}=\rho (C).
\end{align}

Likewise, one has that the received codeword $u+u_0$ is a deep
hole of $C$ if and only if
\begin{align}
\min_{v\in C}\{d(u+u_0, v)\}=\rho (C).
\end{align}

Since
$$
\{d(u+u_0, v)| v\in C\}=\{d(u+u_0, v+u_0)| v\in C\},
$$
it follows that
\begin{align}
\min_{v\in C}\{d(u+u_0, v)\}=\min_{v\in C}\{d(u+u_0, v+u_0)\}.
\end{align}
But $d(u+u_0, v+u_0)=d(u, v)$ for any codeword $u_0$. Hence (2.3)
tells us that
\begin{align}
\min_{v\in C}\{d(u+u_0, v)\}=\min_{v\in C}\{d(u, v)\}.
\end{align}

Now from (2.1), (2.2) and (2.4), one can deduce that
$u$ is a deep hole of $C$ if and only if $u+u_0$ is
a deep hole of $C$ as one desires.
So Lemma 2.1 is proved. $\hfill{\Box}$\\

{\it Remark 2.1.} We should point out that if the
codeword $u_0$ is not in $C$, then Lemma 2.1 is not true.\\

In what follows, we let
$$P_{k-1}:=\{f(x)\mid f(x)\in{\bf F}_q[x],\deg f(x)\leq{k-1}\}.$$
We have the following result. \\

\noindent{\bf Lemma 2.2.} {\it Let $\#(D)=q-l$ and let $u\in{\bf F}_q^{q-l+1}$ and
$v\in{\bf F}_q^{q-l+1}$ be two codewords with $u(x)$ and $v(x)$ being the Lagrange
polynomials of $u$ and $v$. If $u(x)=\lambda v(x)+f_{\le{k-2}}(x)$,
where $\lambda\in{\bf F}_q^*$ and $f_{\le{k-2}}(x)\in {\bf F}_q[x]$ is a polynomial
of degree at most $k-2$, then
$$d(u, {\rm GPRS}_q(D, k)) = d(v, {\rm GPRS}_q(D, k)).$$
Further, $u$ is a deep hole of ${\rm GPRS}_q(D, k)$ if and only if
$v$ is a deep hole of ${\rm GPRS}_q(D, k).$}\\

{\it Proof.} Since $u(x)=\lambda v(x)+f_{\leq k-2}(x)$, we have
$u(D)=\lambda v(D)+f_{\leq k-2}(D)$ and $c_{k-1}(u(x))=\lambda c_{k-1}(v(x))$.
By the definition of Hamming distance, we know that for any code
$C$ over ${\bf F}_q$, if $u$ and $v$ are codewords of $C$, then
$$
d(u, v)=d(u+w, v+w)=d(\lambda u, \lambda v)
$$
hold for any codeword $w$ of $C$ and any $\lambda \in {\bf F}_q^*$.
Then from the definition of error distance and noticing that
$u=(u(D), c_{k-1}(u(x)))$, we can deduce immediately that

\begin{align*}
&d(u,{\rm GPRS}_q(D,k))\\
=&\min_{g\in P_{k-1} } d(u,(g(D),c_{k-1}(g(x))))\\
=&\min_{g\in P_{k-1} } d(({u}(D),c_{k-1}(u(x))),(g(D),c_{k-1}(g(x))))\\
=&\min_{g\in P_{k-1} } d((\lambda v(D)+f_{\leq k-2}
(D), c_{k-1}(u(x))),(g(D),c_{k-1}(g(x))))\\
=&\min_{g\in P_{k-1} } d((\lambda v(D)+f_{\leq k-2}
(D),\lambda c_{k-1}(v(x))),(g(D),c_{k-1}(g(x))))\\
=&\min_{g\in P_{k-1} } d((\lambda v(D)+f_{\leq k-2}
(D),\lambda c_{k-1}(v(x))),(g(D)+f_{\leq k-2}(D),c_{k-1}(g(x))))\\
=&\min_{g\in P_{k-1} } d((\lambda v(D),\lambda c_{k-1}(v(x))),(g(D),c_{k-1}(g(x))))\\
=&\min_{g\in P_{k-1} } d((\lambda v(D),\lambda c_{k-1}(v(x))),(\lambda g(D),\lambda c_{k-1}(g(x))))
({\rm since} \ \lambda\in{\bf F}_q^* )\\
=&\min_{g\in P_{k-1} } d((v(D), c_{k-1}(v(x))),(g(D), c_{k-1}(g(x))))\\
=&d((v(D), c_{k-1}(v(x))), {\rm GPRS}_q(D,k))\\
=&d(v,{\rm GPRS}_q(D,k))
\end{align*}
as required. The proof of Lemma 2.2 is complete. $\hfill{\Box}$\\

For a linear $[n, k]$ code $C$ with $n$ and $k$ being the
length and dimension of $C$, respectively, we define the
{\it generator matrix}, denoted by $G$, to be the $k\times n$
matrix of the form $G:=(g_1,\ldots, g_k)^T$, where
$\{g_1,\ldots, g_k\}$ is a basis of $C$ as a vector space.
Since $D=\{y_1,\ldots, y_{q-l}\}$, the following
$k\times (q-l+1)$ matrix
\begin{align}
\left(
\begin{array}{cc}
 1(D) & c_{k-1}(1)\\
 x(D) & c_{k-1}(x)\\
 \vdots & \vdots \\
 x^{k-2}(D) & c_{k-1}(x^{k-2})\\
 x^{k-1}(D) & c_{k-1}(x^{k-1})
\end{array}
\right)
=\left(
\begin{array}{cccc}
 1 &\ldots &1 & 0 \\
 y_1 &\ldots &y_{q-l} & 0\\
 \vdots & \vdots & \vdots&\vdots \\
 y_1^{k-2} & \ldots &y_{q-l}^{k-2} & 0\\
 y_1^{k-1} & \ldots &y_{q-l}^{k-1} & 1
\end{array}
\right)
\end{align}
forms a generator matrix of ${\rm GPRS}_q(D,k)$.
For the purpose of this paper, we will choose
the above matrix as the generator matrix of ${\rm GPRS}_q(D,k)$. \\

\noindent{\bf Lemma 2.3.} \cite{[V]} {\it Let $C$ be an $[n,k]$ linear code and
$G$ be the generator matrix of $C$. Then $C$ is a MDS code if and only if
any $k$ distinct columns of $G$ are linear independent over finite field ${\bf F}_q$.}\\

Throughout this paper, for any nonempty set
$\{\gamma_1, \ldots, \gamma_n\}\subset{\bf F}_q$,
we define the {\it Vandermonde  determinant}, denoted by $V(\gamma_1, \ldots, \gamma_n)$,
as follows:
$$
V(\gamma_1, \ldots, \gamma_n):=\det\left(
\begin{array}{ccc}
1&\dots&1\\
\gamma_1&\dots&\gamma_n\\
\vdots & \vdots&\vdots \\
\gamma_1^{n-1}&\dots&\gamma_n^{n-1}
\end{array}
\right).$$
We have the following well-known result.\\

\noindent{\bf Lemma 2.4.} \cite{[V]} {\it One has}
$$V(\gamma_1, \ldots, \gamma_n)=\prod\limits_{1\leq i<j\leq n}
(\gamma_j-\gamma_i).$$

In the following, we show that the generalized projective
Reed-Solomon code is a MDS code. \\

\noindent{\bf Lemma 2.5.} {\it Let $D\subset{\bf F}_q$. Then ${\rm GPRS}_q(D,k)$
is a $[q-l+1,k]$ MDS code over finite field ${\bf F}_q$.}\\

{\it Proof.} Let $G$ be the generator matrix of ${\rm GPRS}_q(D,k)$
given in (2.5). Write $G:=(G_1, \ldots, G_{q-l+1})$. Let $i_1,\ldots,i_k$
be arbitrary $k$ distinct integers such that $1\le i_1<\cdots <i_k\leq q-l+1.$
We claim that $\det(G_{i_1}, \ldots, G_{i_k})\neq0$ which will be proved in
what follows.

If $i_k\leq q-l$, then it follows that
$$\det(G_{i_1},\ldots,G_{i_k})=V(y_{i_1},\ldots,y_{i_k})
=\prod\limits_{1\leq t<s\leq k}(y_{i_s}-y_{i_t})\neq 0.$$
The claim is true in this case.

If $i_k=q-l+1$, then by expanding the determinant according to
the last column, we arrive at
$$\det(G_{i_1},\ldots,G_{i_k})=V(y_{i_1},\ldots,y_{i_{k-1}})
=\prod\limits_{1\leq t<s\leq k-1}(y_{i_s}-y_{i_t})\neq0.$$
The claim is proved in this case.

Now by the claim, we can derive that any $k$ columns of
the generator matrix $G$ is linear independent.
Then ${\rm GPRS}_q(D,k)$ is a MDS code by Lemma 2.3.

This concludes the proof of Lemma 2.5. $\hfill{\Box}$\\

The following result about the relation between
the covering radius and minimum distance of ${\rm GPRS}_q(D, k)$
will play a key role in this paper which is due to D\"ur \cite{[D]}.\\

\noindent{\bf Lemma 2.6. }\cite{[D]} {\it Let $D$ be
a proper subset of ${\bf F}_q$. Then}
$$\rho({\rm GPRS}_q(D,k))=d({\rm GPRS}_q(D,k))-1.$$\\

The following result is well known.\\

\noindent{\bf Lemma 2.7. }\cite{[ZLC]}
{\it Let $G$ be a generator matrix of a MDS code $C =[n, k]$ over
the finite field ${\bf F}_q$. If the covering radius $\rho(C)=n-k$,
then a received codeword $u\in{\bf F}_q^n$ is a deep hole of $C$
if and only if the $(k+1)\times n$ matrix
$\left(
\begin{array}{c}
 {G}\\
 {u}
\end{array}
\right)$
can be served as the generator matrix of another MDS code.}\\

In what follows, we show a result on the zero-sum problem
in the finite field of odd characteristic.\\

\noindent{\bf Lemma 2.8.} {\it Let $q=p^s$ with $p$ being an
odd prime number and $k$ be an integer with $2\leq k\leq {q-3}$.
Then there exist a subset $I\subseteq {\bf F}_q^*$ with
$\#(I)=k$ such that $\sum\limits_{z\in I}z=0$.}\\

{\it Proof.} Since $p$ is an odd prime number, it follows that for any
$z\in {\bf F}_q^*$, one has $-z\in {\bf F}_q^*$ and $z\ne -z$ since
$2z\ne 0$. But $|{\bf F}_q^*\setminus \{z, -z\}|=q-3\ge 2$ since
$q\ge k+3\ge 5$. Now one can pick $z'\in {\bf F}_q^*\setminus \{z, -z\}$.
Then $-z'\in {\bf F}_q^*\setminus \{z, -z, z'\}$ since $2z'\ne 0$.
Continuing in this way, we finally arrive at
\begin{align}
{\bf F}_q^*=\{z_1,-z_1,\cdots, z_{\frac{q-1}{2}},-z_{\frac{q-1}{2}}\}.
\end{align}
We consider the following cases.

{\it Case 1.} $2\mid k$. In this case, we let
$I=\{z_1,-z_1,\cdots,z_{\frac{k}{2}},-z_{\frac{k}{2}}\}$.
Then $I\subset{\bf F}_q^*$ and we have
$$\sum_{z\in I}z=\sum_{i=1}^{\frac{k}{2}}(z_i+(-z_i))=0$$
as desired. Lemma 2.8 holds if $2\mid k$.

{\it Case 2.} $2\nmid k$. Then $k\ge 3$ and so $q\ge 7$ since
$2\leq k\leq {q-3}$. We claim that there are three distinct
elements $z', z'', z'''\in {\bf F}_q^*$ such that
$z'+z''+z'''=0$, which will be proved by dividing into
the following three subcases.

{\it Case 2.1.} $p=3$. We pick a $z'\in {\bf F}_q^*$. Then $3z'=0$,
$-z'\ne 0$ and $2z'\ne 0$. The latter implies that $z'\ne -z'$.
Since $p=3$ and $q\ge 7$, we deduce that $q\ge 3^2=9$.
Thus $|{\bf F}_q^*\setminus \{z', -z'\}|=q-3\ge 6$. So we can
choose a $z''\in {\bf F}_q^*\setminus \{z', -z'\}$. But $2z''\ne 0$.
Hence $-z''\in {\bf F}_q^*\setminus \{z', -z', z''\}$. It implies
that $z'+z''\ne 0$, namely, $z'+z''\in {\bf F}_q^*$. Furthermore,
we have that $z'+z''$ is not equal to anyone of the four elements
$z', -z', z''$ and $-z''$. That is,
$z'+z''\in {\bf F}_q^*\setminus \{z', -z', z'', -z''\}$.
Hence $-(z'+z'')\in {\bf F}_q^*\setminus \{z', -z', z'', -z'', z'+z''\}$.
Therefore there are three distinct elements $z', z''$ and $-(z'+z'')$
in ${\bf F}_q^*$ such that their sum equals zero.
The claim holds in this case.

{\it Case 2.2.} $p=5$. Take a $z'\in {\bf F}_q^*$. Then $5z'=0$
and none of $z', 2z', 3z'$ and $4z'$ equals zero. It follows that
the four elements $z', -z', 2z', -2z'$ are pairwise distinct.
Since $q\ge 7>5$, one must have $q\ge 5^2=25$. Thus
$|{\bf F}_q^*\setminus \{z', -z', 2z', -2z'\}|=q-5\ge 20$. So we can
choose $z''\in {\bf F}_q^*\setminus \{z', -z', 2z', -2z'\}$. Then
$-z''\in {\bf F}_q^*\setminus \{z', -z', 2z', -2z'\}$ and $z'+z''\ne 0$.
The latter one tells us that $-(z'+z'')\in {\bf F}_q^*$.
Obviously, $-z''\ne z''$ since $2z''\ne 0$. Hence
$-z''\in {\bf F}_q^*\setminus \{z', -z', 2z', -2z', z''\}$.

Furthermore, we can deduce that $z'+z''$ is not equal to
any of $z', -z', 2z', -2z', z''$ and $-z''$. This infers that
$-(z'+z'')\in {\bf F}_q^*\setminus \{z', -z',2z', -2z', z'', -z'', z'+z''\}$
since $2(z'+z'')\ne 0$.
Therefore we can find three distinct elements $z', z''$ and $-(z'+z'')$
in ${\bf F}_q^*$ such that their sum equals zero.
The claim holds in this case. The claim is proved in this case.

{\it Case 2.3.} $p\ge 7$. Then $le\ne 0$ for any integer $l$
with $1\le l\le 6$, where $e$ stands for the identity of
the group ${\bf F}_q^*$. Since $e\ne 0, 4e\ne 0$ and
$5e\ne 0$, we have $e\ne 2e, e\ne -3e$ and $2e\ne -3e$.
So there are three different elements $e, 2e, -3e$
in ${\bf F}_q^*$ such that their sum is equal to
zero as one desires. The claim is true in this case.

Now by the claim, we know that there are three integers
$i_1, i_2$ and $i_3$ such that $1\le i_1<i_2<i_3\le \frac{q-1}{2}$
and $z_{i_1}+z_{i_2}+z_{i_3}=0$.

If $q=7$, then letting $I=\{z_{i_1}, z_{i_2}, z_{i_3}\}$ gives us the
desired result.

If $q>7$, then
${\bf F}_q^*\setminus\{\pm z_{i_1}, \pm z_{i_2}, \pm z_{i_3}\}$
is nonempty. By (2.6), we obtain that
\begin{align}
\nonumber &{\bf F}_q^*\setminus\{\pm z_{i_1}, \pm z_{i_2}, \pm z_{i_3}\}\\
=&\{\pm z_1, ..., \pm z_{i_1-1}, \pm z_{i_1+1}, ..., \pm z_{i_2-1},
\pm z_{i_2+1}, ..., \pm z_{i_3-1}, \pm z_{i_3+1}, ..., \pm z_{\frac{q-1}{2}}\}.
\end{align}
Since $2\nmid k$, $k-3$ is even.
Evidently, the sum of the first $k-3$ elements on the right hand side
of (2.7) is equal to zero because $z_i+(-z_i)=0$ for all integers
$1\le i\le \frac{q-1}{2}$. Then the first $k-3$ elements on the right
hand side of (2.7) together with the three elements
$z_{i_1}, z_{i_2}, z_{i_3}$ gives us the desired result.
Thus Lemma 2.8 is true if $2\nmid k$.

This completes the proof of Lemma 2.8. $\hfill{\Box}$\\

For any positive integer $x$, we define its $p$-adic valuation,
denoted by $v_p(x)$, to be the largest exponent $r$
such that $p^r$ divides $x$. In the conclusion of this section,
We provide the following characterization on the divisibility
of certain binomial coefficients by the prime number $p$.\\

\noindent{\bf Lemma 2.9.} {\it Let $q$ be a power of the odd prime $p$
and let $t$ be an integer with
$2\le t\le q-1$. Then $v_p\big(\binom{q-2}{t-1}\big)=v_p(t).$
Consequently, the binomial coefficient $\binom{q-2}{t-1}$
is divisible by $p$ if and only if $t$ is a multiple of $p$.}

{\it Proof.}
Clearly, one has
$$\binom{q-2}{t-1}=(q-t)\prod\limits_{i=2}^{t-1}\frac{q-i}{i}.$$
Therefore
\begin{align}
\nonumber v_p(\binom{q-2}{t-1})=&\sum\limits_{i=2}^{t-1}v_p(\frac{q-i}{i})+v_p(q-t)\\
=&\sum\limits_{i=2}^{t-1}(v_p(q-i)-v_p(i))+v_p(q-t).
\end{align}

Since $q$ is a power of $p$, it follows that for any positive
integer $i$ with $i<q$, one has $v_p(i)<v_p(q)$, and so
\begin{align}
v_p(q-i)=v_p(i).
\end{align}
Then from (2.8) and (2.9) one derives that
$$v_p(\binom{q-2}{t-1})=v_p(q-t)=v_p(t)$$
as required. It then follows that
$p\mid\binom{q-2}{t-1}$ if and only if $p\mid t$.
So Lemma 2.9 is proved. $\hfill{\Box}$

\section{\bf Proofs of Theorems 1.4 and 1.6}

In this section, we use the lemmas presented in the previous
to give the proofs of Theorems 1.4 and 1.6. At first,
we show Theorem 1.4. \\

{\it Proof of Theorem 1.4.}
Since $\deg(u(x))=k$, one may let
$u(x)=\lambda x^k+\nu x^{k-1}+f_{\leq{k-2}}(x)$
with $\lambda\in {\bf F}_q^*,\nu\in{\bf F}_q$
and $f_{\leq{k-2}}(x)\in{\bf F}_q[x]$
being a polynomial of degree at most $k-2$.
Then $(u(D), c_{k-1}(u(x)))=(u(D), \nu)$.
By Lemma 2.2, we have that $(u(D), \nu)$
is a deep hole of the generalized projective Reed-Solomon
code ${\rm GPRS}_q(D, k)$ if and only if
$(\lambda ^{-1}u(D), \lambda ^{-1}\nu)$ is a deep hole of ${\rm GPRS}_q(D, k)$.
But
$$\lambda ^{-1}u(x)=w_k(x)+r_k(x),$$
where $w_k(x):=x^k$ and
$$r_k(x):=\lambda ^{-1}\nu x^{k-1}+\lambda ^{-1}f_{\leq{k-2}}(x).$$
Then one has
$$(\lambda ^{-1}u(D), \lambda ^{-1}\nu)
=(w_k(D)+r_k(D), \lambda^{-1}\nu )=(w_k(D), 0)+(r_k(D), \lambda^{-1}\nu).$$

Since $\deg r_k(x)\le k-1$, by the definition of ${\rm GPRS}_q(D, k)$
we have $(r_k(D), \lambda^{-1}\nu)\in {\rm GPRS}_q(D, k)$.
Then it follows from Lemma 2.1 that $(\lambda ^{-1}u(D), \lambda ^{-1}\nu)$
is a deep hole of ${\rm GPRS}_q(D, k)$ if and only if $(w_k(D), 0)$
is a deep hole of ${\rm GPRS}_q(D, k)$. Then we can deduce that
$(u(D), c_{k-1}(u(x)))$ is a deep hole of
${\rm GPRS}_q(D, k)$ if and only if $(w_k(D), 0)$
is a deep hole of ${\rm GPRS}_q(D, k)$.

We denote ${\bar w_k}:=(w_k(D),0)$. Let $G$ be
the generator matrix of ${\rm GPRS}_q(D,k)$
as given in (2.5). Then we have
\begin{align*}
\left(
\begin{array}{c}
 G\\
 {\bar w_k}
\end{array}
\right)
=&\left(
\begin{array}{cccc}
 1 &\ldots &1&0 \\
 y_1 & \ldots &y_{q-l}&0 \\
 \vdots & \vdots & \vdots&\vdots\\
  y_1^{k-2} & \ldots &y_{q-l}^{k-2}&0\\
 y_1^{k-1} & \ldots &y_{q-l}^{k-1}&1\\
 y_1^k & \ldots &y_{q-l}^k & 0
\end{array}
\right)\\
:=&(\bar G_1, \ldots, \bar G_{q-l}, \bar G_{q-l+1}).
\end{align*}

Now we pick $k+1$ distinct integers with
$1\leq j_1<\cdots<j_{k+1}\leq{q-l+1}$.

{\it Case 1.} $j_{k+1}\leq{q-l}$. Then one has
\begin{align*}
\det(\bar G_{j_1}, \ldots, \bar G_{j_{k+1}})
=&\det \left(
\begin{array}{ccc}
 1 &\ldots &1 \\
 y_{j_1} & \ldots &y_{j_{k+1}}\\
 \vdots & \vdots & \vdots\\
 y_{j_1}^{k-1} & \ldots &y_{j_{k+1}}^{k-1}\\
 y_{j_1}^k & \ldots &y_{j_{k+1}}^k
\end{array}
\right)\\
=&V(y_{j_1},\ldots,y_{j_{k+1}})\\
=&\prod\limits_{1\leq t<s\leq {k+1}}(y_{j_s}-y_{j_t})\neq 0.
\end{align*}

{\it Case 2.} $j_{k+1}={q-l+1}$. We can compute and get that
\begin{align}
\nonumber\det(\bar G_{j_1}, \ldots, \bar G_{j_k}, \bar G_{j_{q-l+1}})
=&\det\left(
\begin{array}{cccc}
1 &\ldots &1 & 0\\
 y_{j_1} &\cdots &y_{j_k} & 0\\
 \vdots & \vdots &\vdots &\vdots\\
  y_{j_1}^{k-2} & \ldots &y_{j_k}^{k-2}&0\\
 y_{j_1}^{k-1} & \ldots &y_{j_k}^{k-1}& 1\\
   y_1^k& \ldots &y_{j_k}^k& 0
\end{array}
\right)\\
=&-\det\left(
\begin{array}{ccc}
1 &\ldots &1\\
 y_{j_1} &\cdots &y_{j_k}\\
 \vdots & \vdots &\vdots \\
 y_{j_1}^{k-2} & \ldots &y_{j_k}^{k-2}\\
y_{j_1}^{k} & \ldots &y_{j_k}^{k}\\
\end{array}
\right).
\end{align}

Now we introduce an auxiliary polynomial $g(y)$
as follows:

$$g(y)=\det\left(
\begin{array}{cccc}
1 &\ldots &1 & 1\\
 y_{j_1} &\cdots &y_{j_k} & y\\
 \vdots & \vdots &\vdots &\vdots\\
 y_{j_1}^{k-1} & \ldots &y_{j_k}^{k-1}& y^{k-1}\\
  y_{j_1}^{k} &\ldots &y_{j_k}^{k}& y^k
\end{array}
\right).
$$
Then Lemma 2.4 tells us that
$$g(y)=\Big(\prod\limits_{1\leq s<t\leq k}(y_{j_t}-y_{j_s})\Big)
\prod\limits_{i=1}^{k}(y-y_{j_i}):=\sum\limits_{i=0}^{k}a_iy^i.
$$
It infers that
\begin{align}
a_{k-1}=-\Big(\sum\limits_{i=1}^ky_{j_i}\Big)
\prod\limits_{1\leq s<t\leq k}(y_{j_t}-y_{j_s}).
\end{align}
But
\begin{align}
a_{k-1}=-\det\left(
\begin{array}{ccc}
1 &\ldots &1\\
 y_{j_1} &\cdots &y_{j_k}\\
 \vdots & \vdots &\vdots \\
 y_{j_1}^{k-2} & \ldots &y_{j_k}^{k-2}\\
y_{j_1}^{k} & \ldots &y_{j_k}^{k}\\
\end{array}
\right).
\end{align}
Finally, (3.1) together with (3.2) and (3.3) gives us that
\begin{align}
\det(\bar G_{j_1}, \ldots, \bar G_{j_k}, \bar G_{j_{q-l+1}})
=-\Big(\sum\limits_{i=1}^ky_{j_i}\Big)\prod\limits_{1\leq s<t\leq k}(y_{j_t}-y_{j_s}).
\end{align}

By Lemma 2.5, we know that ${\rm GPRS}_q(D,k)$
is a $[q-l+1, k]$ MDS code which implies that
$$d({\rm GPRS}_q(D,k))=q-l+1-k+1=q-l-k+2.$$
Then by Lemma 2.6, one can deduce that
\begin{align}
\rho({\rm GPRS}_q(D,k))=&d({\rm GPRS}_q(D,k))-1\\
\nonumber=&q-l-k+2-1\\
\nonumber=&q-l+1-k.
\end{align}
It then follows immediately from Lemma 2.7 that
$\bar w_k=(w_k(D), 0)$ is a deep hole of the generalized
projective Reed-Solomon code ${\rm GPRS}_q(D, k)$
if and only if the
$(k+1)\times (q-l+1)$ matrix
$\left(
\begin{array}{c}
 G\\
 {\bar w_k}
\end{array}
\right)$
can be served as the generator matrix of a MDS code,
if and only if any $k+1$ columns of
$\left(
\begin{array}{c}
 G\\
 {\bar w_k}
\end{array}
\right)$ are linear independent, if and only if for any
$1\le j_1<\cdots<j_{k+1}\le q-l+1$, one has
\begin{align}
\det(\bar G_{j_1}, \ldots, \bar G_{j_{k+1}})\neq 0.
\end{align}
By the discussion in Cases 1 and 2, (3.4) tells us that (3.6) holds
if and only if for any $1\le j_1<\cdots<j_{k}\le q-l$,
one has $\sum\limits_{i=1}^ky_{j_i}\neq 0$.
Hence we can derive that $(w_k(D), 0)$ is a deep
hole of the generalized projective Reed-Solomon code
${\rm GPRS}_q(D, k)$ if and only if the sum
$\sum\limits_{y\in I}y$ is nonzero for any subset
$I\subseteq D$ with $\#(I)=k$ as desired.

Finally, we can conclude that $(u(D), c_{k-1}(u(x)))$
is a deep hole of the generalized projective Reed-Solomon
code ${\rm GPRS}_q(D, k)$ if and only if the sum
$\sum\limits_{y\in I}y$ is nonzero for any subset
$I\subseteq D$ with $\#(I)=k$.

This finishes the proof of Theorem 1.4. $\hfill{\Box}$\\

We can now use Theorem 1.4 to show Theorem 1.6.\\

{\it Proof of Theorem 1.6.} Let $l=1$ and $a_1=0$.
Then $D={\bf F}_q^*$. By Lemma 2.8, there exist a subset
$I\subseteq {\bf F}_q^*$ with $\#(I)=k$
such that $\sum\limits_{y\in I}y=0$. It then follows
from Theorem 1.4 that
the received word $(u({\bf F}_q^*),c_{k-1}(u(x)))$=$(u({\bf F}_q^*),\gamma)$
is not a deep hole of the primitive projective Reed-Solomon
code ${\rm PPRS}_q({\bf F}_q^*, k)$. Therefore Theorem 1.6 is proved. $\hfill{\Box}$

\section{\bf Proofs of Theorems 1.5 and 1.7}

In this section, we give the proof of Theorems 1.5 and 1.7.
We begin with the proof of Theorem 1.5. \\

{\it Proof of Theorem 1.5.}
First of all, we note that $j$ is an integer with
$1\leq j\leq l$. We introduce a polynomial $f_j(x)$ as follows:
$$f_j(x)=(x-a_j)^{q-2},$$
and define a codeword ${\bar f_j}$ associated to $f_j(x)$ by
$${\bar f_j}:=(f_j(D), c_{k-1}(f_j(x))).$$
Then $u_j(x)=\lambda_j f_j(x)+\nu_jx^{k-1}+f_{\le k-2}^{(j)}(x)$
which implies that
\begin{align}
c_{k-1}(u_j(x))=\lambda_j c_{k-1}(f_j(x))+\nu_j.
\end{align}
It follows from (4.1) that
\begin{align*}
&(u_j(D), c_{k-1}(u_j(x)))\\
=&(\lambda_j f_j(D)+\nu_jx^{k-1}(D)+f_{\le k-2}^{(j)}(D), \lambda_jc_{k-1}(f_j(x))+\nu_j)\\
=&(\lambda_j f_j(D), \lambda_jc_{k-1}(f_j(x)))+(\nu_jx^{k-1}(D)+f_{\le k-2}^{(j)}(D), \nu_j)\\
=&\lambda_j\bar f_j+(\nu_jx^{k-1}(D)+f_{\le k-2}^{(j)}(D), \nu_j).
\end{align*}
But
$$\deg (\nu_jx^{k-1}(x)+f_{\le k-2}^{(j)}(x))\le k-1$$
and
$$c_{k-1}(\nu_jx^{k-1}(x)+f_{\le k-2}^{(j)}(x))=\nu_j.$$
Hence
$$(\nu_jx^{k-1}(D)+f_{\le k-2}^{(j)}(D),\nu_j)\in {\rm GPRS}_q(D, k).$$
It follows from Lemmas 2.1 and 2.2 that the received word
$(u_j(D), c_{k-1}(u_j(x)))$ is a deep hole of the generalized
projective Reed-Solomon code ${\rm GPRS}_q(D, k)$ if and only if
$\bar f_j$ is a deep hole of ${\rm GPRS}_q(D, k)$.

Let $G$ be the generator matrix of ${\rm GPRS}_q(D,k)$
as given in (2.5). Since $y_{i}\ne a_{j}$ for all integers
$i$ with $1\leq i\leq q-l$, we have
$(y_{i}-a_{j})^{q-2}=(y_{i}-a_{j})^{-1}$. It then follows that
\begin{align}
\left(
\begin{array}{c}
G\\
{\bar f_j}
\end{array}
\right)
=&\left(
\begin{array}{cccc}
 1 &\ldots &1 & 0 \\
 y_1 &\ldots &y_{q-l} & 0\\
 \vdots & \vdots & \vdots&\vdots \\
 y_1^{k-2} & \ldots &y_{q-l}^{k-2} & 0\\
 y_1^{k-1} & \ldots &y_{q-l}^{k-1} & 1\\
 (y_1-a_j)^{-1} &\ldots &(y_{q-l}-a_j)^{-1} & c_{k-1}(f_j)
\end{array}
\right)\\
\nonumber:=&(\hat G_1, \ldots, \hat G_{q-l+1}).
\end{align}

On the other hand, from Lemma 2.7 we can deduce that
$\bar f_j=(f_j(D), c_{k-1}(f_j(x)))$ is a deep hole
of the generalized projective Reed-Solomon code
${\rm GPRS}_q(D, k)$, if and only if
$\left(
\begin{array}{c}
 G\\
 {\bar f_j}
\end{array}
\right)$
generates a MDS code, by Lemma 2.3, if and only if
any $k+1$ columns of
$\left(
\begin{array}{c}
 G\\
 {\bar f_j}
\end{array}
\right)$ are linear independent, if and only if for all
$k+1$ integers $j_1, \ldots, j_{k+1}$ with
$1\le j_1<\cdots<j_{k+1}\le q-l+1$, one has
\begin{align}
\det(\hat G_{j_1}, \ldots, \hat G_{j_{k+1}})\neq 0.
\end{align}

In what follows, we choose arbitrarily $k+1$
integers $j_1,\ldots,j_{k+1}$ such that
$1\le j_1<\cdots <j_{k+1}\leq q-l+1$.
Consider the following two cases.

{\it Case 1.} $j_{k+1}\ne q-l+1$. Then
$k+1\le j_{k+1}\le q-l$ and by (4.2), one has
\begin{equation*}
(\hat G_{j_1}, \ldots, \hat G_{j_{k+1}})
=\left(
\begin{array}{ccc}
 1 &\ldots &1 \\
 y_{j_1} &\ldots &y_{j_{k+1}}\\
 \vdots &\vdots &\vdots\\
 y_{j_1}^{k-1} &\ldots &y_{j_{k+1}}^{k-1}\\
(y_{j_1}-a_j)^{-1} &\ldots &(y_{j_{k+1}}-a_j)^{-1}
\end{array}
\right).
\end{equation*}
Thus one can deduce that
\begin{align*}
&\det(\hat G_{j_1}, \ldots, \hat G_{j_{k+1}})\\
=&\Bigg(\prod\limits_{i=1}^{k+1}(y_{j_i}-a_j)^{-1}\Bigg)
\det\left(
\begin{array}{ccc}
y_{j_1}-a_j &\ldots &y_{j_{k+1}}-a_j\\
 y_{j_1}(y_{j_1}-a_j)&\ldots &y_{j_{k+1}}(y_{j_{k+1}}-a_j)\\
 \vdots & \vdots  &\vdots\\
 y_{j_1}^{k-1}(y_{j_1}-a_j)&\ldots &y_{j_{k+1}}^{k-1}(y_{j_{k+1}}-a_j)\\
 1 & \ldots &1 \\
\end{array}
\right)\\
=&\Bigg(\prod\limits_{i=1}^{k+1}(y_{j_i}-a_j)^{-1}\Bigg)
\det\left(
\begin{array}{ccc}
y_{j_1} &\ldots &y_{j_{k+1}}\\
 \vdots & \vdots  &\vdots\\
 y_{j_1}^{k}&\ldots &y_{j_{k+1}}^{k}\\
 1 & \ldots &1 \\
\end{array}
\right)\nonumber
\end{align*}
\begin{align*}
=&(-1)^k\Bigg(\prod\limits_{i=1}^{k+1}({y_{j_i}-a_{j}})^{-1}\Bigg)
V(y_{j_1},\ldots, y_{j_{k+1}})\nonumber\\
=&(-1)^{k}\Bigg(\prod\limits_{i=1}^{k+1}({y_{j_i}-a_{j}})^{-1}\Bigg)
\prod\limits_{1\leq s<t\leq{k+1}}(y_{j_t}-y_{j_s})\ne 0
\end{align*}
since $y_{j_1}, \ldots, y_{j_{k+1}}$ are pairwise distinct.

{\it Case 2.} $j_{k+1}=q-l+1$. Then $1\le j_1<\cdots <j_{k}\leq q-l$.
From (4.2) and Lemma 2.4, we can deduce that
\begin{align}
\nonumber&\det(\hat G_{j_1},\ldots,\hat G_{j_{k}},\hat G_{j_{q-l+1}})\\
\nonumber=&\det\left(
\begin{array}{cccc}
1 &\ldots &1 & 0\\
 y_{j_1} &\cdots &y_{j_k} & 0\\
 \vdots & \vdots &\vdots &\vdots\\
  y_{j_1}^{k-2} & \ldots &y_{j_k}^{k-2}& 0\\
 y_{j_1}^{k-1} & \ldots &y_{j_k}^{k-1}& 1\\
  (y_{j_1}-a_j)^{-1} &\ldots &(y_{j_k}-a_j)^{-1}& c_{k-1}(f_j(x))
\end{array}
\right)\\
\nonumber=& c_{k-1}(f_j(x)) \det\left(
\begin{array}{ccc}
1 &\ldots &1\\
 y_{j_1} &\cdots &y_{j_k}\\
 \vdots & \vdots &\vdots \\
 y_{j_1}^{k-1} & \ldots &y_{j_k}^{k-1}\\
\end{array}
\right)-\det\left(
\begin{array}{ccc}
1 &\ldots &1\\
 y_{j_1} &\cdots &y_{j_k}\\
 \vdots & \vdots &\vdots \\
 y_{j_1}^{k-2} & \ldots &y_{j_k}^{k-2}\\
(y_{j_1}-a_j)^{-1} &\ldots &(y_{j_k}-a_j)^{-1}
\end{array}
\right)\\
\nonumber=&c_{k-1}(f_j(x))V(y_{j_1},\ldots, y_{j_{k}})-\det\left(
\begin{array}{ccc}
1 &\ldots &1\\
 y_{j_1} &\cdots &y_{j_k}\\
 \vdots & \vdots &\vdots \\
 y_{j_1}^{k-2} & \ldots &y_{j_k}^{k-2}\\
(y_{j_1}-a_j)^{-1} &\ldots &(y_{j_k}-a_j)^{-1}
\end{array}
\right)\\
\nonumber=&c_{k-1}(f_j(x)) V(y_{j_1},\ldots, y_{j_{k}})
-\Bigg(\prod\limits_{i=1}^{k}(y_{j_i}-a_j)^{-1}\Bigg)
\det\left(
\begin{array}{ccc}
y_{j_1} &\ldots &y_{j_{k}}\\
 \vdots & \vdots  &\vdots\\
 y_{j_1}^{k-1}&\ldots &y_{j_{k}}^{k-1}\\
 1 & \ldots &1 \\
\end{array}
\right)\\
\nonumber=&c_{k-1}(f_j(x)) V(y_{j_1},\ldots, y_{j_{k}})
+(-1)^k \Bigg(\prod\limits_{i=1}^{k}(y_{j_i}-a_j)^{-1}\Bigg)
V(y_{j_1},\ldots, y_{j_{k}})\\
\nonumber=&\Big(c_{k-1}(f_j(x))+(-1)^k\prod\limits_{i=1}^{k}(y_{j_i}-a_j)^{-1})\Big)
\prod\limits_{1\leq s<t\leq k}(y_{j_t}-y_{j_s})\\
=&\Big(c_{k-1}(f_j(x))+\prod\limits_{i=1}^{k}(a_j-y_{j_i})^{-1})\Big)
\prod\limits_{1\leq s<t\leq k}(y_{j_t}-y_{j_s}).
\end{align}

Now from Cases 1 and 2, we can deduce by (4.4) that (4.3)
holds for all $k+1$ integers $j_1, \ldots, j_{k+1}$ with
$1\le j_1<\cdots<j_{k+1}\le q-l+1$
if and only if for all integers $j_1, \ldots, j_k$ with
$1\le j_1<\cdots<j_{k}\le q-l$, one has
$$c_{k-1}(f_j(x))+\prod\limits_{i=1}^{k}(a_j-y_{j_i})^{-1}\neq 0,$$
which is equivalent to
\begin{align}
c_{k-1}(f_j(x))\prod\limits_{i=1}^{k}(a_j-y_{j_i})+e\neq 0.
\end{align}

Since $f_j(x)=(x-a_j)^{q-2}$, the binomial theorem gives us that
$$c_{k-1}(f_j(x))=\binom {q-2}{k-1}(-a_j)^{q-k-1}.$$
Then one derives that (4.5) holds for all integers
$j_1, \ldots, j_k$ with $1\le j_1<\cdots<j_{k}\le q-l$
if and only if the following is true:
\begin{align}
\binom {q-2}{k-1}(-a_j)^{q-k-1}\prod\limits_{i=1}^{k}(a_j-y_{j_i})+e\neq 0,
\end{align}
or equivalently,
\begin{align}
\binom {q-2}{k-1}a_j^{q-k-1}\prod\limits_{i=1}^{k}(y_{j_i}-a_j)+e\neq 0
\end{align}
since $q$ is odd. In other words, $\bar f_j=(f_j(D), c_{k-1}(f_j(x)))$
is a deep hole of the generalized projective Reed-Solomon code
${\rm GPRS}_q(D, k)$ if and only if the sum
$$\binom{q-2}{k-1}a_j^{q-1-k}\prod\limits_{y\in I}(y-a_j)+e$$
is nonzero for any subset $I\subseteq D$ with $\#(I)=k$.
Hence the desired result follows immediately.
The first part is proved.

Now we show the second part. Let $k\equiv 0\pmod p$.
Then by Lemma 2.9, we have
$\binom{q-2}{k-1}\equiv 0\pmod p$. So one can write
$\binom{q-2}{k-1}=p\Delta$ with $\Delta $ being a positive integer.
Then
$$c_{k-1}(f_j(x))=\binom{q-2}{k-1}(-a_j)^{q-k-1}=p\Delta (-a_j)^{q-k-1}=0.$$
It then follows that
$$\binom{q-2}{k-1}a_j^{q-1-k}\prod\limits_{y\in I}(y-a_j)+e=e\ne 0$$
for any subset $I\subseteq D$ with $\#(I)=k$. So it follows from
the first part that the received codeword $(u_j(D), c_{k-1}(u_j(x)))$
is a deep hole of ${\rm GPRS}_q(D, k)$ if $k\equiv 0\pmod p$.
The second part is proved.

The proof of Theorem 1.5 is complete. $\hfill{\Box}$\\

We can now present the proof of Theorem 1.7
as the conclusion of this paper.\\

{\it Proof of Theorem 1.7.} Letting $l=1$ and $a_1=0$
gives us that $D={\bf F}_q^*$.

If $p|k$, then by Theorem 1.5, one knows that
the received word $(u({\bf F}_q^*),\delta)$
is a deep hole of the primitive projective Reed-Solomon
code ${\rm PPRS}_q({\bf F}_q^*, k)$.

If $p\nmid k$, then it follows from $a_1=0$ that
$$\binom{q-2}{k-1}a_j^{q-1-k}\prod\limits_{y\in I}(y-a_j)+e
=0\cdot \binom{q-2}{k-1}\prod\limits_{y\in I}(y-a_j)+e=e\ne 0$$
for any subset $I\subseteq D$ with $\#(I)=k$. Hence
$(u({\bf F}_q^*),\delta)$ is a deep hole of ${\rm PPRS}_q({\bf F}_q^*, k)$.

This completes the proof of Theorem 1.7. $\hfill{\Box}$

\end{document}